%% file: manuscript_ogiermann_et_al_clean.tex
\documentclass[
        parskip = half, %
        twoside,
        footinclude = false, %
        ]{scrartcl}
\pdfoutput=1
\usepackage{ifthen} %

\newboolean{boldeqnitalic}
\setboolean{boldeqnitalic}{true}
\usepackage{stylefiles/defdb1} %

\usepackage{stylefiles/rublkm_paper} %

\usepackage{todonotes}

\usepackage{algorithm}
\usepackage{algpseudocode}

\usepackage{siunitx}
\DeclareSIUnit{\molar}{M}

\usepackage{pgfplots}
\pgfplotsset{compat=1.17}
\usepgfplotslibrary{groupplots}
\usepackage{tikzscale}

\usepackage{multirow}

\makeatletter
\newcommand*{\authormark}{}
\newcommand*{\markauthor}[1]{%
  \renewcommand{\authormark}{#1}%
  \ignorespaces
}
\newcommand*{\titlemark}{}
\newcommand*{\marktitle}[1]{%
  \renewcommand{\titlemark}{#1}%
  \ignorespaces
}
\makeatother

\markauthor{D. Ogiermann, D. Balzani,  L.E. Perotti}
\marktitle{An Explicit Local Space-Time Adaptive Framework for Monodomain Models}

\clearpairofpagestyles
\rohead{\pagemark}
\rehead{\pagemark}

\lohead{\titlemark}
\lehead{\authormark}

\KOMAoptions{
   headsepline = true
}
\setkomafont{pagehead}{\normalfont\normalcolor\small}

\newcommand{\average}[1]{\ensuremath{\left\{#1\right\}} }
\newcommand{\jump}[1]{\ensuremath{[\![#1]\!]} }

\usepackage[normalem]{ulem}

\usepackage{mathtools}
\DeclarePairedDelimiter\ceil{\lceil}{\rceil}

\DeclarePairedDelimiter\set{\lbrace}{\rbrace}

\begin{document}

\input{sec_01_clean}
\section*{Acknowledgement}

We thank the MFEM team for making their framework publicly available and providing quick help on technical details. 
Daniel Balzani acknowledges financial funding from the Deutsche Forschungsgemeinschaft (DFG), project ID~465228106. 
Furthermore, the free and open source software community, as well as scientific discussions with Maximilian Köhler and Hendrik Dorn are greatly acknowledged.

\bibliographystyle{stylefiles/references_style}
\bibliography{references_clean}

\end{document}

%% file: sec_01_clean.tex
\begin{center}

{\LARGE An Explicit Local Space-Time Adaptive Framework for Monodomain Models in Cardiac Electrophysiology}

\vspace{5mm}

Dennis Ogiermann$^{1,\star}$, Daniel Balzani$^{1}$, Luigi E. Perotti$^{2}$

\vspace{3mm}

{\small $^1$Chair of Continuum Mechanics, Ruhr University Bochum,\\ 
Universitätsstraße~150, 44801~Bochum, Germany}\\[3mm]

{\small $^2$Department of Mechanical and Aerospace Engineering, University of Central Florida,\\
Orlando, USA}\\[3mm]

\vspace{3mm}

{\small ${}^{\star}$E-mail address of corresponding author: dennis.ogiermann@rub.de}

\vspace{10mm}

\begin{minipage}{15.0cm}
\textbf{Abstract}\hspace{3mm}
We present a new explicit local space-time adaptive framework to decrease the time required for  monodomain simulations for cardiac electrophysiology. 
Based on the localized structure of the steep activation wavefront in solutions to monodomain problems, the proposed framework adopts small time steps and a tree-based adaptive mesh refinement scheme only in the regions necessary to resolve these localized structures.
The time step and mesh adaptation selection process is fully controlled by a combination of local error indicators.
The main contributions of this work consist in the introduction of a primal symmetric interior penalty formulation of the monodomain model and an efficient algorithmic strategy to manage local time stepping for its temporal discretization.
In a first serial implementation of this framework, we report decreases in wall-clock time between 2 and 20 times with respect to an optimized implementation of a commonly used numerical scheme, showing that this framework is a promising candidate to accelerate monodomain simulations of cardiac electrophysiology.
\end{minipage}
\end{center}

\medskip{}
\textbf{Keywords:} Synchronous Local Time-Stepping, Adaptive Mesh Refinement, Space-Time Adaptivity, Discontinuous Galerkin, Cardiac Electrophysiology, Computational Cardiology

\section{Introduction}

Computational cardiology has gained wide attention during the last decade as it can provide valuable insights in understanding cardiac function and dysfunction.
For example, numerical simulations have been used to investigate drugs in virtual screenings~\cite{PasBriLuRohHerGalGreBueRod:2017:hsd,PasZhoTroBriBueRod:2021:vas}, to understand the heart activation sequence and origin of the electrocardiogram~\cite{KriPerBorAjiFriPonWeiQuKluEnnGar:2014:smv,PerKriBorEnnKlu:2015:rsv,CarBueMinZemRod:2016:hva,SanDAMafCaiPriKraAurPot:2018:sav,OgiBalPer:2021:ema,CamLawDroMinWanGraBurRod:2021:iva,MosWulLewHorPerStrMenKraOdeSee:2022:cmr}, and to assist in diagnosis~\cite{PerPonKriBalEnnKlu:2017:mui,GraDobHei:2019:cmw,PerVerMouCorLoeBalEnn:2021:ecs}, prognosis~\cite{AroAliTra:2019:rpa}, and therapy planning~\cite{AroAliTra:2019:rpa,GraDobHei:2019:cmw,HubHeiZonvandeTorHosDelvan:2018:wnm}.
For the wide adoption of cardiac models in these use-cases, two aspects are very important.
First, the adopted numerical scheme has to be verified, the model has to be properly validated, and uncertainties should be quantified.
In this way the model predictions can be trustworthy~\cite{Joh:2016:SCV,PatGra:2018:vtm,OgiBalPer:2021:ema}.
Second, the simulations should be reasonably fast, especially if the simulation models are to be deployed in a clinical setting for near-real time predictions. 
Strict requirements on time step length and grid size have been shown to be the major bottleneck in many simulations (c.f., e.g.,~\cite{NieKerBenBerBerBraCheClaFenGarHeiLanMalPatPlaRodRoySacSeeSkaSmi:2011:vct,KriSarKlu:2013:nqo,PezHakSun:2015:sea} for some analyses) stemming from the characteristic fast, steep traveling wavefronts.

In light of these limitations on simulation efficiency, a significant effort has been devoted to increase the speed of cardiac electrophysiology models.
For simulations requiring only an approximation of the activation sequence, a simplified model based on the Eikonal wave approximation, namely the Reaction-Eikonal model, has been developed~\cite{NeiCamPraNieBisVigPla:2017:ece}.
This approach has been shown to be helpful in rapid ECG evaluations~\cite{PezKalPotPriAurKra:2017:ERA} when coupled with a proper lead field approach~\cite{OgiBalPer:2021:ema,NagEspGilGseSanPlaDosLoe:2023:cpm}.
Although this approximation is the fastest method available for many forward simulations, it is of limited use in several scenarios, e.g. when investigating the initiation and dynamics of fibrillation, drug interaction, and mechano-electrical feedback or defibrillation.
In such scenarios, bidomain~\cite{Tun:1978:bmd} and monodomain models are the most common choice.
We will focus on the latter in this paper. %

When the solution of the monodomain (or bidomain) models is required, several strategies have been proposed. 
\cite{ColPav:2004:psr,FraDeuErdLanPav:2006:AST,MunPav:2009:dsa,WonGokKuh:2013:cmc,HurHen:2014:gfv,BarHuyPavSca:2022:pns} have investigated monolithic schemes.
While monolithic schemes appear as an obvious choice, they do not take direct advantage of the problem specific structure to optimize the solution strategy, and this results in sub-optimal performance (see, e.g.,~\cite{LinGerJahLoeWeiWie:2023:ets}).
Since their introduction, a significant amount of work has instead focused on the development and optimization of operator splitting schemes~\cite{QuGar:1999:aas,HeiGasFerRod:2010:csa,HeiFerDobRod:2010:amfa,CriMasRav:2013:fam,CerSpi:2018:hom,GreSpi:2019:gis,GomOliLobdos:2020:amm,MouPue:2021:dae,OgiPerBal:2023:sea} as they have been shown to be efficient and simple to implement~\cite{KriSarKlu:2013:nqo,SacMarBurMeiConWeb:2018:peg}.

Recently~\citet{KabCheFen:2019:ris} have shown that explicit time discretization schemes for cardiac electrophysiology can outperform implicit and also operator splitting based schemes.
An interesting takeaway is that the CFL condition from the explicit discretization of the diffusion operator and the time step length to resolve the cell model at the wavefront do not differ significantly for sufficiently fine grids to resolve the steep wavefront.
Similar observations on the CFL condition have been made in the dual-adaptive explicit scheme by~\citet{MouPue:2021:dae}.
In addition, in our previous work~\cite{OgiPerBal:2023:sea}, we have studied how the time step length can be efficiently and automatically adapted during different phases of the simulation, providing evidence on the temporal locality of the problem.
However, all these schemes still use a small time step length everywhere in the domain during depolarization. 
Hence, they key question in this paper consists in how to construct a numerical method which utilizes the spatial and temporal locality of monodomain problems and which is at least competitive with respect to existing methods. 

\citet{KraKra:2016:elt} approached this question by constructing a numerical method using space-time elements in combination with a block-adaptive mesh.
Interestingly, in their work, \citet{KraKra:2016:elt} report that their implementation is unable to outperform operator splitting schemes, likely due to the overhead associated with assembling and solving the space-time problem and multiple grid adaptions per time step.
We start from~\citet{KraKra:2016:elt}'s idea of using the local nature of the monodomain problem, but instead of using an implicit scheme for the time discretization, we develop an improved space-time adaptive scheme with local explicit time stepping.
Further, we take inspiration from local time stepping schemes found in computational fluid dynamics~\cite{Kri:2010:elt,GasStaHinAtaMun:2015:sta} and existing local adaptive time stepping techniques found in efficient operator splitting schemes for cardiac electrophysiology problem~\cite{QuGar:1999:aas,Whi:2007:PDA,OgiPerBal:2023:sea}.
Regarding the motivation for spatial adaptivity in the proposed scheme, the efficiency of using adaptive mesh refinement algorithms to exploit the spatial locality has already been demonstrated to increase the efficiency of numerical simulations of monodomain problems~\cite{CheGreHen:2000:sam,CheGreHen:2003:Est,BelForBou:2009:anm,HeiFerDobRod:2010:amfa,ArtBisKay:2013:esc,DicKraKraPot:2014:dal,SacMarBurMeiConWeb:2018:peg,HoeBerKroPfaChaWal:2018:ahd}.

The rest of this work is structured as follows.
In Section~\ref{sec:methods}, we begin by stating the monodomain problem in its strong form.
Then, we derive the weak primal form of the monodomain problem in the symmetric interior penalty framework to transform the system of partial differential equations into a system of ordinary differential equation (Section~\ref{sec:spatial-discretization}).
This constitutes the first novel contribution of this work.
Subsequently, before presenting the time marching scheme for this system of ordinary differential equations, we describe how our framework adapts the mesh using a tree-based refinement to control the spatial errors (Section~\ref{sec:spatial-adaptivity}).
In Section~\ref{sec:s-lts} the second novel contribution of this work, i.e., the adaptive ODE solver for the semi-discrete monodomain problem, is described and analyzed.
The description of the method is followed by the conduction velocity benchmark proposed by~\citet{NieKerBenBerBerBraCheClaFenGarHeiLanMalPatPlaRodRoySacSeeSkaSmi:2011:vct} to verify our scheme (Section~\ref{sec:niederer-benchmark}).
Further, an additional test based on computing the full revolution of a spiral wave is proposed in Section~\ref{sec:spiral-wave-benchmark}.
This test aims to verify that the wavefront follows its actual trajectory instead of the spatial discretization, a phenomenon known as lattice pinning (e.g., see~\cite[Fig.~39]{FenCheHasEva:2002:mms}).
Our final benchmark (Section~\ref{sec:lv-benchmark}) consists in an idealized left ventricular geometry without microstructure to investigate the performance of the proposed scheme on more complex geometries (than rectangular domains).
We conclude with a discussion of limitations, future directions, and applications for the proposed framework.

\section{Methods}
\label{sec:methods}

The monodomain model is a common framework to describe the spatiotemporal dynamics of tissue electrophysiology at the continuum level.
This approach considers the evolution of two fields: the transmembrane potential field $\varphi$ and a state field $\boldsymbol{s}$ representing the state of the cellular model.
The monodomain model defined on a domain $\Omega$ can be written as

\begin{subequations}
    \begin{align}
        \chi  C_{\textrm{m}} \partial_t \varphi &= \nabla \cdot \boldsymbol{\kappa} \nabla \varphi - \chi I'(\varphi, \boldsymbol{s}, t) & \textrm{in} \: \Omega \, ,  \label{eq:monodomain_A} \\
        \partial_t \boldsymbol{s} &= \mathbf{g}(\varphi, \boldsymbol{s}) & \textrm{in}  \: \Omega \, , \label{eq:monodomain_B} \\
        0 &= \boldsymbol{\kappa} \nabla \varphi \cdot \mathbf{n} & \mathrm{on} \: \partial \Omega \, ,
        \label{eq:monodomain_C}
    \end{align}
    \label{eq:monodomain}
\end{subequations}

together with admissible initial conditions and a cellular ionic model to determine $I'$ and $\mathbf{g}$. 
In eq.~\eqref{eq:monodomain}, ${\boldsymbol{\kappa}}$ denotes the conductivity tensor, $\varphi_{\mathrm{m}}$ is the transmembrane potential field, $\chi$ is the volume to membrane surface ratio, $C_{\mathrm{m}}$ is the membrane capacitance, and $I'(\varphi, \boldsymbol{s}, t) := I_{\textrm{ion}}(\varphi, \boldsymbol{s}) + I_{\textrm{stim}}(t)$ denotes the sum of the ionic current due to the cell model and the applied stimulus current, respectively.
For the numerical solution of this set of differential equations, a new adaptive scheme is proposed in this section with the following main components: (i) A spatial discretization based on discontinuous Galerkin, which allows to decouple the time-stepping on each element; (ii) A spatial error indicator enabling spatial adaptivity; and (iii) An efficient adaptive local time-stepping scheme. 
These components are described in detail in the following. 

\subsection{Spatial Discretization}
\label{sec:spatial-discretization}

We first discretize the monodomain model (eq.~\eqref{eq:monodomain}) in space with a discontinuous Galerkin scheme.
The key idea motivating the use of a discontinuous Galerkin (DG) discretization is that we can easily decouple the time derivatives in each element allowing different time step lengths per element and thus, local time-stepping.
First, we approximate the domain $\Omega$ with $N_{\mathrm{E}}$ discrete, non-overlapping elements such that
\[
    \Omega \approx \bigcup^{N_{\mathrm{E}}}_{e=1} \Omega_e \, .
\]
Next, we state the monodomain model's weak form in DG primal formulation. 
We start by multiplying eq.~\eqref{eq:monodomain_A} with a test function~$\delta\varphi$ and integrating over the domain $\Omega$.
This yields the following standard equation, also used as starting point for the continuous Galerkin approach
\begin{equation}
    \int_{\Omega} \partial_t \varphi \delta\varphi \,\mathrm{d}V = \int_{\Omega} \nabla \cdot \boldsymbol{D} \nabla \varphi \delta\varphi - I(\varphi, \boldsymbol{s}, t) \delta\varphi \,\mathrm{d}V \qquad\forall \delta\varphi \in \mathbb{T}\, 
    \label{eq:contGal}
\end{equation}
where $\mathbb{T}$ is a suitable test space, $\boldsymbol{D} = \boldsymbol{\kappa} / \chi C_{\textrm{m}}$ is the diffusion tensor, and $I = I'/C_{\textrm{m}}$. 
For the discontinuous Galerkin approach, element-wise discontinuous shape functions are considered leading to independent degrees of freedom (DOFs) at the element boundaries. 
After partial integration of Equation~\eqref{eq:contGal} and use of appropriate integral theorems, we obtain the following integral equation
\begin{equation}
    \begin{aligned}
        \int_{\Omega} \partial_t \varphi \delta\varphi \,\mathrm{d}V &= -\int_{\Omega} \boldsymbol{D} \nabla \varphi \cdot \nabla \delta\varphi + I(\varphi, \boldsymbol{s}, t) \delta\varphi \,\mathrm{d}V & \\ 
        & \quad + \int_{\Gamma} \average{\vec{\sigma}}\cdot\jump{\delta\varphi} + \jump{ \varphi} \average{\boldsymbol{D} \nabla \delta\varphi} \mathrm{d}S & \forall \delta\varphi \in \mathbb{T} \, , \\
    \end{aligned}
    \label{eq:monodomain-primal-formulation}
\end{equation}
cf.~\citet{ArnBreCocMar:2002:uad}. 
Here, $\Gamma$ denotes the union of the inter-element faces (not including the outer boundary of the domain), $\vec{\sigma}$ is the \textit{numerical flux}, $\average{\bullet} = (\bullet^+ + \bullet^-)/2$ denotes the \textit{average}, and $\jump{\bullet} = \bullet^+ \boldsymbol{n}^+ + \bullet^- \boldsymbol{n}^-$ is the \textit{jump} operator with $\bullet^+$ and $\bullet^-$ denoting the quantities on the two sides of a face and $\boldsymbol{n}$ the outward pointing unit normal. 
For the solution field to be physically meaningful, the jump of the transmembrane potential should vanish in the limit as the mesh is refined. 
One approach to accomplish this consists in setting the numerical flux to 
\begin{equation}
    \vec{\sigma} = \average{\boldsymbol{D} \nabla \varphi} + \gamma \average{\frac{\boldsymbol{n}\cdot\boldsymbol{D}\boldsymbol{n}}{h}} \jump{\varphi} \, ,
    \label{eq:monodomain-numerical-flux}
\end{equation}
where $h$ is the characteristic size of the element and $\gamma$ is the penalty parameter leading to vanishing jumps of the solution field. 
In doing so, we recover an anisotropic version of the \textit{symmetric interior penalty Galerkin} scheme~\cite{DouDup:1976:ipp}.
We finalize the spatial semi-discretization by choosing piecewise discontinuous polynomials for ansatz and test space, where we denote the discretized solution field at the nodes by the vectors $\boldsymbol{\tilde{\varphi}}$ and $\tilde{\boldsymbol{s}}$.
The resulting system of ODEs can then be written as
\begin{equation}
    \begin{aligned}
        \boldsymbol{M} \mathrm{d}_t \boldsymbol{\tilde{\varphi}} &= \boldsymbol{K}\boldsymbol{\tilde{\varphi}} + \mathbf{N}(\boldsymbol{\tilde{\varphi}}, \boldsymbol{\tilde{s}}) \, , %
        \label{eq:monodomain-ode-massmatrix-form}
    \end{aligned}
\end{equation}
where $\boldsymbol{M}$ is the block-diagonal mass matrix (one dense block per element), $\boldsymbol{K}$ is the blocked diffusion matrix (one block per element and per shared face with neighboring elements), and $\mathbf{N}(\boldsymbol{\tilde{\varphi}}, \boldsymbol{\tilde{s}})$ is the nonlinear block-form associated to the ionic part of the cell model.
It should be noted that, as it is true for $\boldsymbol{M}$, also the blocks of $\mathbf{N}$ do not couple across elements, as all involved quantities and the test functions are element-local.
Since the mass matrix is block-diagonal we can efficiently invert it for approximation spaces with a low number of DOFs, or when applying diagonalization techniques, resulting in the final form of the spatial discretization
\begin{equation}
    \begin{aligned}
        \mathrm{d}_t \boldsymbol{\tilde{\varphi}} &= \boldsymbol{M}^{-1} \boldsymbol{K}\boldsymbol{\tilde{\varphi}} + \boldsymbol{M}^{-1} \mathbf{N}(\boldsymbol{\tilde{\varphi}}, \boldsymbol{\tilde{s}}) \, , \\
        \mathrm{d}_t \boldsymbol{\tilde{s}} &= \boldsymbol{g}(\boldsymbol{\tilde{\varphi}}, \boldsymbol{\tilde{s}}) \, .
        \label{eq:monodomain-ode-form}
    \end{aligned}
\end{equation}
If line, quadrilateral, and hexahedral elements are used to discretize the spatial domain $\Omega$, we select the orthogonal Lagrange polynomials with nodes based on the roots of Legendre polynomials for the ansatz and test spaces, such that the nodes are in the interior of the elements.
The integrals in the matrices and the nonlinear form are approximated based on Gauss-Legendre quadrature.
This choice for the integrals and the function spaces diagonalizes the mass matrix and decouples the nonlinear form point-wise.
In principle, these choices allow to utilize higher order approximations for the spatial approximation efficiently. 
However, it is also possible to use other element types -- possibly at the cost of a dense inverse mass matrix per element.
In this context, preliminary numerical testing of nodal positions and quadrature rules other than the ones listed above for line, quadrilateral, and hexahedral elements resulted in a slightly reduced efficiency. 

To implement the described spatial discretization, we use the popular finite element framework MFEM~\cite{AndAndBarBraCamCerDobDudFisKolPazStoTomAkkDahMedZam:2021:mmf} in version 4.5 for management of the mesh (including adaptive mesh refinement described below) and finite element operators.

\subsection{Spatial Adaptivity}
\label{sec:spatial-adaptivity}

\begin{figure}
    \begin{minipage}{0.5\textwidth}
        \centering
        \input{figures/amr.tikz}
    \end{minipage}
    \begin{minipage}{0.5\textwidth}
        \includegraphics[width=\textwidth]{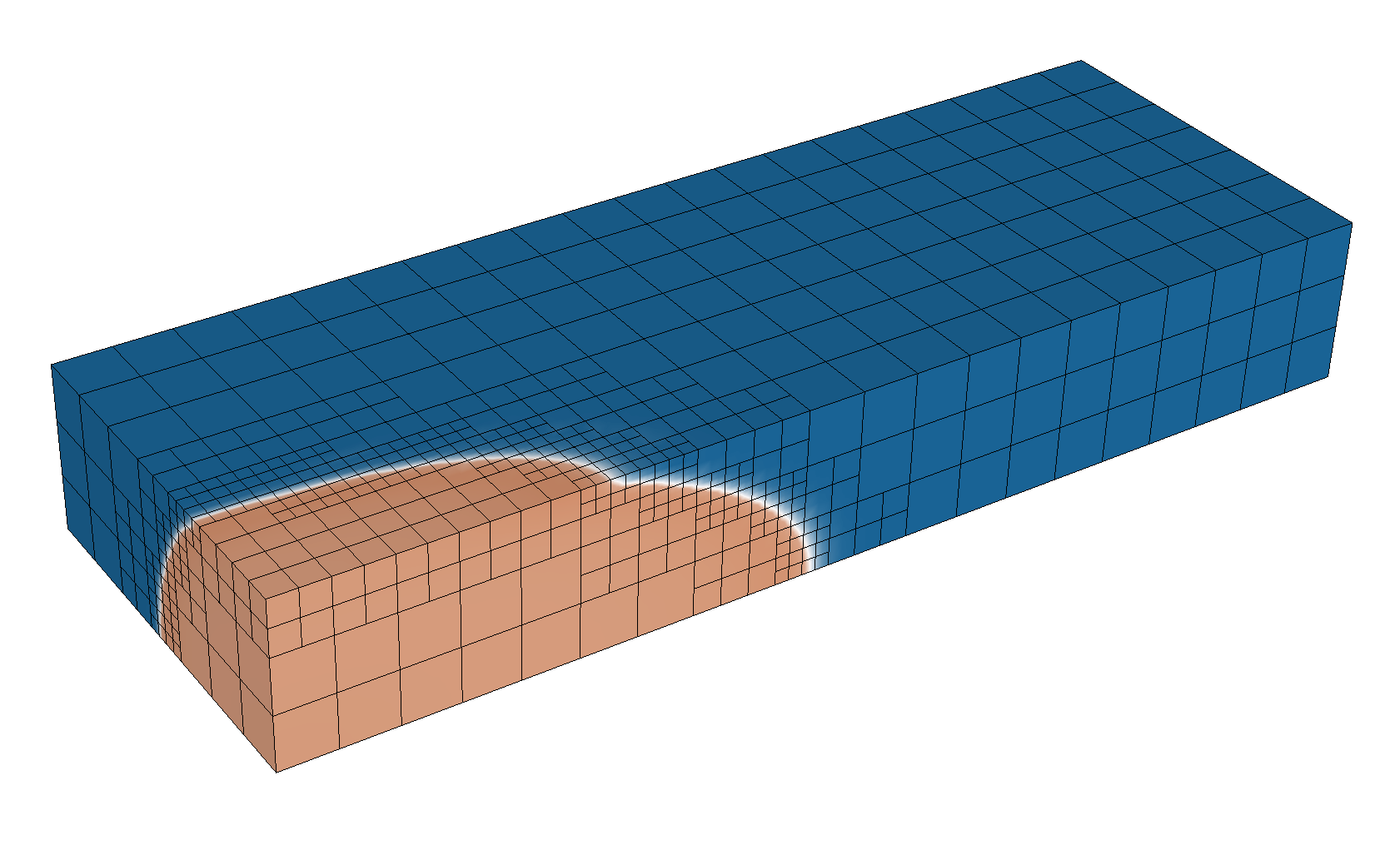}
    \end{minipage}
        \caption{
        \textbf{Left}: Example of tree-based adaptive mesh refinement on a single quadrilateral element.
        In this case, refinement operations add four children (i.e., reading the figure left to right) while coarsening operations collapse them into their parent elements (i.e., reading the figure right to left).
        \textbf{Right}: Adaptive mesh refinement with hexahedral elements during the conduction velocity benchmark~\cite{NieKerBenBerBerBraCheClaFenGarHeiLanMalPatPlaRodRoySacSeeSkaSmi:2011:vct} (rendered with GLVis).
        It can be clearly seen that the mesh is properly refined around the wavefront and coarsened everywhere else.
    }
    \label{fig:adaptive-grid}
\end{figure}

Since our temporal discretization strategy will depend on details of the spatial adaptivity, we start by describing the latter.
Using adaptive mesh refinement will allow us to drastically reduce the number of DOFs in each timestep since the wavefront is usually significantly smaller than the full domain size in cardiac tissue and organ simulations, and it is highly localized in space.
The basic strategy consists in classical h-refinement, utilizing a forest of trees to track the non-conforming refinement and coarsening of individual elements.
Refinement operations replace the parent elements with $n$ child elements of the same geometry such that the parent element is partitioned into $n$ equally sized elements, while coarsening reverts this procedure (Fig.~\ref{fig:adaptive-grid}).
For more details on this strategy we refer to the original publication describing this implementation~\cite{CerDobKol:2019:nmr}.
This adaptive mesh strategy is very general and hence not optimized for our specific use-case.
This may lead to suboptimal performance and could be improved in future specific applications.
For the remainder of this work we call the \textit{initial mesh} the \textit{root mesh} and the elements of the initial mesh the \textit{root elements} to reflect on the forest of trees structure of this approach.

Due to the chosen DG approximation, the jump in the flux across neighboring elements is the dominant error component. Indeed the jump in flux can be interpreted as the ``local non-smoothness'' of the piece-wise discontinuous approximation of the smooth wavefront. Based  on this consideration, we estimate the spatial error of the DG approximation with the following Kelly-type error indicator (where the original estimator has been derived in~\cite{KelDeZieBab:1983:pea}) utilizing a penalized jump in the flux over the faces
\begin{equation}
    \eta^{\mathrm{s}}_e = \sqrt{\sum_{F \in \partial\Omega_e} \frac{h_F}{2p} ||W_F\jump{\boldsymbol{D} \nabla \varphi \cdot \boldsymbol{n}}||^2_{L_2(F)}} \, ,
    \label{eq:dg-error-estimate}
\end{equation}
where $h_F$ denotes the characteristic size of the face, $p$ is the order of the ansatz space, and $W_F$ is the transformation weight of the face (from the reference isoparametric domain to the spatial domain $\Omega_e$), which acts as a penalty factor depending on the face size $h_F$.
Hence, $W_F$ can be interpreted as a helper to prevent over-refinement at the wavefront as it decreases with smaller element size.
The use of $W_F$ is the primary modification over the original Kelly error indicator.
An advantage of this penalty approach over using a limiter for the refinement depth is that it also slightly refines repolarization regions.
This feature could help, for example, to capture early afterdepolarizations, although we have not tested this aspect in our simulations.
The proposed error indicator can also be modified to increase the mesh resolution of the wavefront in slow directions by dropping the conductivity tensor from the estimate. %

For the adaptive mesh refinement process, we introduce two threshold parameters $\tau_{\textrm{refine}}$ and $\tau_{\textrm{coarsen}}$, where an element is refined if its estimated error is below the threshold $\tau_{\textrm{refine}}$ and it is coarsened if all children of a parent element are above the threshold $\tau_{\textrm{coarsen}}$.

\subsection{Synchronous Local Time Stepping}
\label{sec:s-lts}

Our main goal is to develop a \textit{synchronous local time stepping} (S-LTS) scheme for the monodomain model, which is based on the scheme described in~\cite{Kri:2010:elt}.
The notion of ``synchronous'' refers to the fact that the time step length on each element is a unit fraction of the global time step. 
In principle, it would be possible to derive a fully asynchronous scheme as presented in~\cite{GasJan:2016:wrs}.
However, this would lead to a significantly more complex implementation and the requirement of appropriate predictor and corrector steps. 
To decouple the time evolution of each element we discretize~\eqref{eq:monodomain-ode-form} in time with a forward Euler scheme.
Also in this regard, it would be possible to utilize a higher order (explicit or exponential) time discretization scheme at the cost of using a predictor-corrector scheme, as for example presented and analyzed in~\cite{Kri:2010:elt,GasJan:2016:wrs} for compressible fluid flow simulations. 
The internal variables are discretized in time with a partitioned first order exponential integrator, which is well-known to be stable in practice~(see, e.g.,\cite{PerVen:2009:egr}).
In this way we recover an analogue of the well-known Rush-Larsen scheme~\cite{RusLar:1978:PAS}.

Let us consider specific time intervals $[t_n, t_{n+1}]$ with $\Delta t_n = t_{n+1} - t_{n}$, such that $t_{n}$ is the time instance where the adaptive refinement is executed. 
We call $\Delta t_n$ the \textit{barrier time step length}, because $t_{n+1}$ (and $t_{n}$) can be seen as the synchronization points where the solution will be guaranteed to be synchronized, and hence it builds a barrier between consecutive time integrations.
For the discretization of diffusion problems using explicit schemes, the Courant-Friedrichs-Lewy (CFL) condition~\cite{CouFriLew:1928:upd} states that the time step length is, up to a scheme-dependent constant, restricted linearly by the reciprocal of the maximum eigenvalue of the conductivity tensor (or physically speaking by the fastest conduction direction) and quadratically by the element size.

As stated by~\citet{RoyBouPie:2020:atm}, the analysis of Rush-Larsen discretization is still problematic.
Hence, we also cannot provide a full formal derivation of the critical time step length for such a mixed scheme, and we do not close this gap in this work.
Instead, we adopt a heuristic approach based on the assumption that the diffusion operator is the limiting factor for the critical time step length.
In our analysis, we start by freezing the solution $\tilde{\boldsymbol{u}}$ for all degrees of freedom associated with an element $e$ and treat it as constant.
We defined the block vector $\boldsymbol{u}(\boldsymbol{x}, t) := [\varphi(\boldsymbol{x}, t), \boldsymbol{s}^\mathrm{T}(\boldsymbol{x}, t)]^\mathrm{T}$ and dropped the arguments to simplify the notation.
We also ignore the nonlinear contributions.
With this information, we investigate the time step length on a single element based on the diffusion part of the stiffness matrix, which we define as:
\begin{equation}
    \boldsymbol{L}_e := \boldsymbol{M}_e^{-1} \boldsymbol{K}_e \, .
\end{equation}
Now, we simply use the eigenvalues of this matrix as a heuristic to restrict the time step length of a given element $e$.
The eigenvalues of this matrix are estimated by Gershgorin discs~\cite{Ger:1931:uae} instead of being computed explicitly.
Accordingly, we obtain the restriction
\begin{equation}
    \Delta t^e 
    \leq \underbrace{\min_{i \in \text{rows of }\boldsymbol{L}_e} \left( \frac{1}{L_{e_{ii}} + \sum_{i \neq j} |L_{e_{ij}}|} \right)}_{:=\mathrm{CFL}_e}
    \label{eq:cfl-monodomain}
\end{equation}
and we refer to~\cite{MouPue:2021:dae} for a detailed analysis of the relation between the spectrum of $\boldsymbol{L}_e$ and the critical time step length of the forward Euler scheme for diffusion problems that motivated this heuristic.
In the presented benchmarks, this estimated restriction on $\Delta t^e$ always resulted in stable time steps.

Next, we focus on the proposed time marching strategy. 
Given the barrier time step length $\Delta t$, we assign a number of substeps $S_e$ to an element.
As discussed above, the CFL condition is the first constraining factor for selecting the time step length of an element.
Hence, a number of substeps can be computed by assigning a time step length satisfying the CFL condition of an element $e$ via the smallest integer $b^{\mathrm{CFL}}_e$ such that
\begin{equation}
    \Delta t \cdot 2^{-b^{\mathrm{CFL}}_e} \leq \mathrm{CFL}_e
    \implies b^{\mathrm{CFL}}_e = \ceil{\log_2(\mathrm{CFL}_e^{-1} \Delta t)} \, , 
\end{equation}
where $\ceil{\bullet}$ denotes the rounding up of $\bullet$ to the next larger integer. 
Assigning powers of 2 to each element forces the scheme to be synchronous: the time step lengths of adjacent elements are guaranteed to be multiples of each other.
This is motivated by the observation that element refinement approximately divides the stable time step length by a factor of four. 

Another important consideration in choosing the element time step length is that cells can undergo localized fast transients in the internal variables (some Markov chain models are prominent examples, see e.g.~\cite{SpiDea:2010:sac} for an in-depth analysis).
These fast transients may be badly resolved when using the time step length based only on the CFL condition, which only states that stability is guaranteed without assuring small approximation errors. 
To reduce the approximation error, we leverage a recurring idea in cardiac electrophysiology: the time step length is adapted as a function of an estimate of the time discretization error of the cell model.
We refer to~\cite{QuGar:1999:aas} and \cite{OgiPerBal:2023:sea} for examples, where the latter is based on the error analysis in~\cite{SpiTor:2016:osb}.
However, this analysis is not directly applicable to the proposed scheme as originally it has been carried out in an operator splitting framework, while our work presents a monolithic framework.
Instead, we will use a modified version of the bound given in the a posteriori error analysis of a monolithic monodomain discretization given by~\cite{RatVer:2019:pee}, which we will call the Ratti-Veroni temporal error indicator or, in short, RV-T indicator.

The RV-T indicator is given by the Bochner-norms (see, e.g., reference~\cite{Eva:2010:pde}, ch. 5.9.2) on an element $e$ in the time interval $[t_a, t_b]$ as
\begin{align}
    \eta^{\mathrm{t}}_{e,I} &= \frac{1}{t_b-t_a}||I(\hat{\tilde{\boldsymbol{u}}}(\boldsymbol{x},t), t) - I(\hat{\tilde{\boldsymbol{u}}}(\boldsymbol{x},t_a), t_a)||_{L_{2}([t_a, t_b], \Omega_e)} \, , \\
    \eta^{\mathrm{t}}_{e,g} &= \frac{1}{t_b-t_a}||g(\hat{\tilde{\boldsymbol{u}}}(\boldsymbol{x},t)) - g(\hat{\tilde{\boldsymbol{u}}}(\boldsymbol{x},t_a))||_{L_{2}([t_a, t_b], \Omega_e)} \, , \\
    \eta^{\mathrm{t}}_e &= \sqrt{{\eta^{\mathrm{t}}}^2_{e,I} + {\eta^{\mathrm{t}}}^2_{e,g} } \, ,
    \label{eq:ratti-veroni-temporal}
\end{align}
where the inner vector norm is simply the Euclidean norm and $\hat{\tilde{\boldsymbol{u}}}$ denotes the time-discretization of the spatially discretized degrees of freedom~$\tilde{\boldsymbol{u}}$. 
We compute the number of time substeps on the element to resolve the local temporal dynamics of the cell as a function of the RV-T indicator $\eta^{\mathrm{t}}_e$, where the time integrals are approximated with the midpoint rule.
We want to highlight the fact that we do not need full evaluations of $I$ and $g$.
It is sufficient to evaluate the fast ionic currents and the components of the right hand side of $g$ in the variables with possibly fast transients.
However, choosing the components to be evaluated (or even just an indicator of these) is model dependent and not explored in this work.
It mainly represents an additional opportunity to further improve efficiency. 
For the remainder of this work we simply use a threshold to compute $b^{\mathrm{cell}}_{e}$ as either $0$ or a fixed number, dependent on the barrier time step length. 
This allows us to compute the number of substeps on each element as 
\begin{equation}
    S_e = 2^{\max(b^{\mathrm{CFL}}_{e}, b^{\mathrm{cell}}_{e})} \, 
    \quad\mbox{with}\quad
    b^{\mathrm{cell}}_{e} = \begin{cases}  
    \ceil{\log_2\left(\Delta t/\overline{\Delta t}\right)} &\mbox{if}\; \eta^{\mathrm{t}}_e > \tau_{\mathrm{cell}}\\
    0 & \mbox{else.}
    \end{cases}
    \label{eq:num-substeps}
\end{equation}
Herein, $\overline{\Delta t}$ and $\tau_{\mathrm{cell}}$ are predefined threshold values. 
This definition ensures stability of the scheme and a sufficient temporal resolution for accuracy at the element level.\\ 

\begin{figure}
    \centering
    \includegraphics[width=\textwidth]{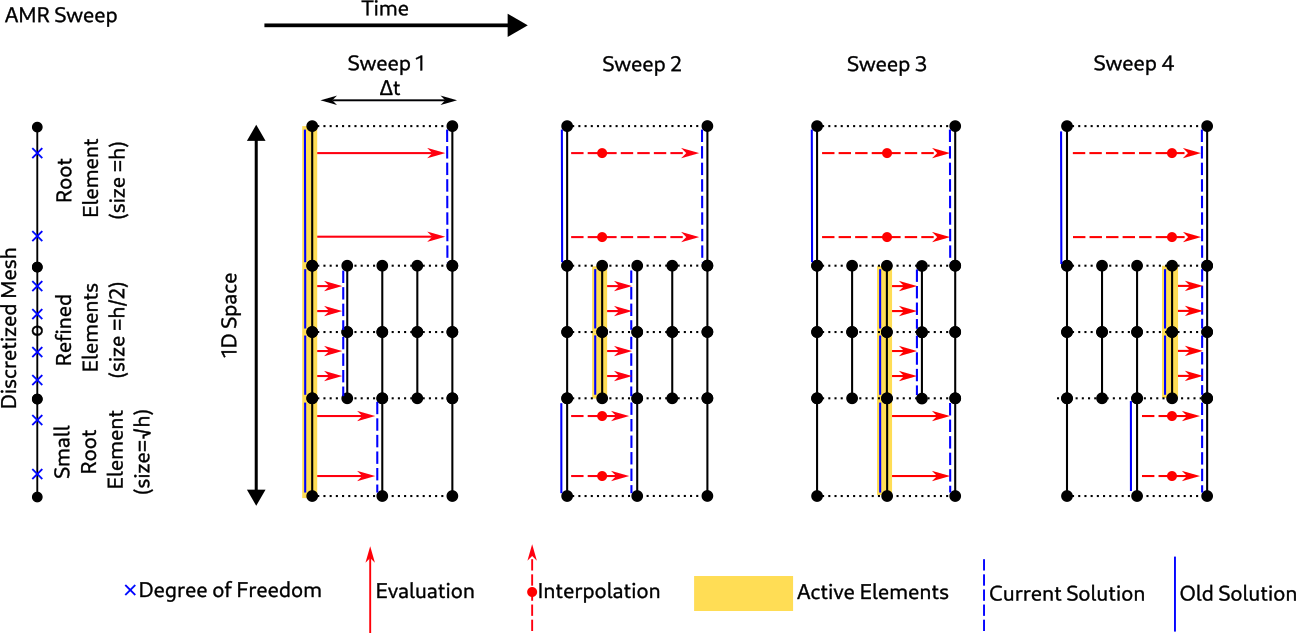}
    \caption{
        Schematic for one barrier time step on a mesh with three root elements.
        Here the center element is refined at the beginning of the time step.
        Lines highlighted in yellow illustrate the active elements where the solution has been computed in the last time substep. 
        Solid red arrows indicate the advance in time whereas dashed red arrows correspond to computing the interpolated values (red circular markers). 
    }
    \label{fig:s-lts-visualized}
\end{figure}

Next we describe the time marching algorithm on an example before describing its efficient implementation in the next paragraph and Algorithm~\ref{alg:multiqueue-s-lts}.
The algorithm is visualized in Fig.~\ref{fig:s-lts-visualized} for a one-dimensional mesh with 4 elements.
Before advancing in time, the mesh is adapted to the current solution with the strategy outlined in Section~\ref{sec:spatial-adaptivity} and the number of substeps per element is computed as described previously in this section. 
The element-wise update of the solution approximation is executed in consecutive sweeps similar to the algorithm described by~\citet{Kri:2010:elt}. 
For these sweeps we use the ordered set of occurring substeps $S_1 \le ... \le S_{N_\text{E}}$ and associated elements to naturally group the elements. 
Note that per construction all substeps $S_e$ are powers of 2 (see Equation~\eqref{eq:num-substeps}), such that there always exists a compatible number of substeps to advance between consecutive groups.
The first sweep advances all elements according to their assigned time step length.
At this time, the solution on each element is advanced to at least $t+\Delta t/S_{N_E}$, marking the current global position in time.
This concludes the first sweep. 
For the next sweep we only consider the elements whose solution is exactly at $t+\Delta t/S_{N_E}$.
These elements are advanced according to their assigned time step length. 
To evaluate the flux integrals across neighboring elements, the solution in the coarser elements at the current sub-timestep is approximated by linear interpolation in time. 
This concludes the second sweep.
The remaining 2 sweeps proceed in similar manner, concluding the time marching from $t_n$ to $t_{n+1}$. 

Regarding an efficient serial implementation of the proposed scheme, we opt to assemble all linear element-wise operators (i.e., $\boldsymbol{K}_e$ and $\boldsymbol{M}^{-1}_e$) before initiating stepping in time.
We note that, as this strategy increases memory usage, it may be semi-optimal if the proposed algorithm was implemented in parallel. 
For the time marching, we need to track three sets of information: 1) The current solution approximation per element; 2) The most recent solution approximation, necessary for time interpolation on the coarser elements; and 3) The update order of the elements.
To track the solution approximations we simply need two vectors, while the element updating is tracked by defining one queue $Q_i$ per substep.
The implementation is reflected in Algorithm~\ref{alg:multiqueue-s-lts}.
Note that it is technically possible to exploit the structure of the algorithm to trivially parallelize the scheme: as the elements' evaluations are independent, we can parallelize over the element queue for each sweep.

\begin{algorithm}
    \caption{Adaptive Multi-Queue Synchronous Local Time Stepping}
    \label{alg:multiqueue-s-lts}
    \begin{algorithmic}[1] %
        \Procedure{A-MQ-S-LTS}{$...$}
            \State Compute spatial error estimates $\eta^{\mathrm{s}}_e$ for each element $e$ via Equation~\eqref{eq:dg-error-estimate}
            \State Refine elements $e$ with $\eta^{\mathrm{s}}_e \geq \tau_{\mathrm{refine}}$ and interpolate solution from old to new mesh
            \State Coarsen elements $e$ with $\eta^{\mathrm{s}}_e \leq \tau_{\mathrm{coarsen}}$ and interpolate solution from old to new mesh
            \State Compute number of substeps $S_e$ via CFL and temporal error estimates $\eta^{\mathrm{t}}_e$
            \State Fill first update queue $Q_0$ with all elements $e$
            \For{$ i \gets 0$ to $[\max_{e}(S_{e})-1]$}
                \For{$ e \in Q_i$}
                    \State Buffer the current solution values for $e$ \Comment{Needed for time interpolation}
                \EndFor
                \For{$ e \in Q_i$}
                    \State Linear interpolate solution of neighboring elements of $e$ to current time
                    \State Update current solution for element $e$ with one Rush-Larsen step
                    \State $j \gets i+\max_{e'}(S_{e'})/S_e$ \Comment{Next substep index}
                    \If{$j \leq \max_{e}(S_{e})-1$}
                        \State Queue element $e$ in $Q_{j}$
                    \EndIf
                \EndFor
            \EndFor
        \EndProcedure
    \end{algorithmic}
\end{algorithm}

\section{Results}

In this section we provide several benchmark calculations using the proposed method to demonstrate its efficiency. 
The studies are performed on an AMD EPYC 7351P and all measurements are carried out with Caliper~\cite{BoeGamBecBreGimLeGPeaSch:2016:cpi} excluding IO operation timings.
All codes were compiled in clang version 15 with optimization level 3 (-O3), native tuning (-mtune=native) and architecture (-march=native), as well as with polyhedral optimizer enabled~\cite{GroGroLen:2012:ppp} (-polly -polly-vectorizer=stripmine). 
Neither a fine tuning of the compiler flags nor profile-guided optimization has been carried out. 
No shared memory or distributed memory parallelization has been implemented for our scheme (i.e., we rely only on the instruction-level parallelism introduced by the compiler). 
We show that the proposed approach outperforms optimized implementations of standard operator splitting schemes~\cite{QuGar:1999:aas}, namely the Lie-Trotter-Godunov scheme with full mass matrix~(e.g.~\cite{KriSarKlu:2013:nqo,OgiPerBal:2023:sea}), by a factor of $\approx 10$. We note that a significant part of the compute time is spent in the AMR loop ($20-50\%$) and performing reassembly of the face and element matrices after refinement and coarsening ($20-30\%)$. These aspects suggest that  further optimizations are possible and we include a preliminary discussion of several possibilities in Sec.~\ref{sec:limitations}.

\subsection{Conduction velocity verification benchmark}
\label{sec:niederer-benchmark}

The benchmark proposed by Niederer et al.~\cite{NieKerBenBerBerBraCheClaFenGarHeiLanMalPatPlaRodRoySacSeeSkaSmi:2011:vct} is a classical verification test for numerical schemes solving the monodomain equations.
Here, instead of the original ten-Tusscher-Panfilov model~\cite{TenPan:2006:asb} used in ~\cite{NieKerBenBerBerBraCheClaFenGarHeiLanMalPatPlaRodRoySacSeeSkaSmi:2011:vct}, we adopt the more recent Pathmanathan-Cordeiro-Gray (PCG) model~\cite{PatCorGra:2019:cuq}, which is a cardiomyocyte electrophysiology model fitted fully against canine data.
The PCG model has been designed for validation and verification purposes and is easier to implement than the ten-Tusscher-Panfilov model, which may help simplifying the reproducibility of the presented results.
Note that the utilization of the PCG model does not significantly affect the solution.
This is in agreement with previous results~\cite{OgiPerBal:2021:alf} comparing the solutions obtained with the classical ten-Tusscher-Panfilov model and the PCG model.
According to this verification benchmark, the analyzed virtual test specimen is a slab with dimensions $\SI{20}{\milli\meter}\times\SI{7}{\milli\meter}\times\SI{3}{\milli\meter}$.
The wave propagation is simulated for $\SI{50}{\milli\second}$ and the diffusion tensor is  $\boldsymbol{D}=\mbox{diag}[0.1334,0.0176,0.0176]$$\si{\milli\meter\squared\per\milli\second}$, 
 where the long axis of the geometry is aligned with the primary eigenvector (associated with the highest diffusivity) of the diffusion tensor.
To obtain a reference solution, the problem is discretized with a Lie-Trotter-Godunov scheme in time and with (conforming) finite elements in space.
Here the time step length is fixed to $\Delta t=\SI{0.01}{\milli\second}$ together with a hexahedral mesh with a characteristic element size of $\SI{0.125}{\milli\meter}$ and a linear Lagrange ansatz, resulting in a total of 229425 DOFs, which is in line with the high quality solution in the original benchmark paper~\cite{NieKerBenBerBerBraCheClaFenGarHeiLanMalPatPlaRodRoySacSeeSkaSmi:2011:vct}. 
For numerical quadrature we used the order 2 Gauss-Legendre rule. 
The linear solver is a preconditioned conjugate gradient scheme~\cite{HesSti:1952:mcg} with Gauss-Seidel preconditioner from hypre~\cite{FalYan:2002:hlh}. 
The solver's relative and absolute tolerances are $10^{-7}$ and $10^{-8}$, respectively, resulting in four to eight iterations per time step.
We confirmed that lower tolerances and a finer grid did not significantly affect the solution.\\ 

\begin{figure}
    \centering
    \begin{minipage}{.29\textwidth}
        \centering
        \includegraphics[width=0.95\textwidth]{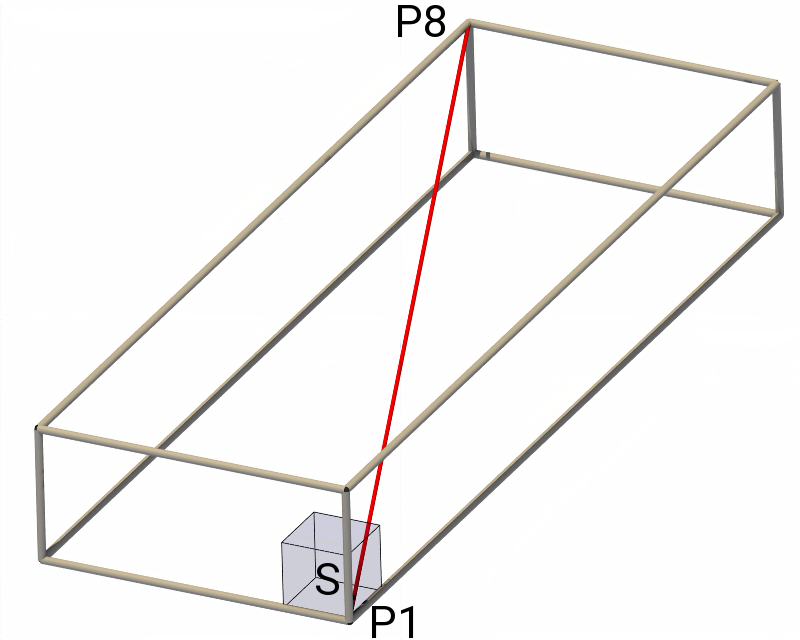}
    \end{minipage}
    \begin{minipage}{.30\textwidth}
        \centering
        \includegraphics[width=0.98\textwidth]{plots/niederer-benchmark-elements.tikz}
    \end{minipage}
    \begin{minipage}{.38\textwidth}
        \centering
        \includegraphics[width=0.98\textwidth]{plots/niederer-analysis.tikz}
    \end{minipage}
    \caption{\textbf{Left:} Schematic illustration of the benchmark setup.
    The stimulus is applied in the region marked with `S' and the red diagonal represents the line along which local activation times are measured.
    \textbf{Center:} The number of time steps for all elements in the mesh over time for the newly proposed S-LTS scheme.
    The root mesh contains 420 elements and the element number peaks at 25159 elements at \SI{27}{\milli\second}.
    This corresponds to the point in time when the wavefront size peaks.
    Hence it can be observed that the AMR follows the wavefront closely since the number of element evaluations directly correlates with the number of elements through the enforcement of the CFL condition.
    \textbf{Right:} Local activation times along the measurement line obtained with the LTG (reference) and S-LTS solutions.
    Note that the reference solution is in good quantitative agreement with the solution reported in the original conduction velocity benchmark study of~\citet{NieKerBenBerBerBraCheClaFenGarHeiLanMalPatPlaRodRoySacSeeSkaSmi:2011:vct}.
    We computed the solution for \SI{50}{\milli\second}.
    The proposed S-LTS method required \SI{81}{\second} while the reference simulation, using an optimized implementation of the classical operator splitting scheme from~\cite{KriSarKlu:2013:nqo,QuGar:1999:aas} (see text for more details), required \SI{1237}{\second} ($\approx \SI{21}{\minute}$), resulting in a speed up of $\approx 15$.
    IO operations as well as the time for the LAT computation were excluded from these time measurements.
    }
    \label{fig:niederer-benchmark}
\end{figure}

For our new S-LTS scheme we use a barrier time step length $\Delta t = \SI{0.15}{\milli\second}$ (see Section~\ref{sec:s-lts}) and a hexahedral root mesh containing 420 elements, corresponding to a characteristic element size of $\SI{1.0}{\milli\meter}$. 
As the threshold for mesh refinement (see Equation~\eqref{eq:dg-error-estimate} in Section~\ref{sec:spatial-adaptivity}) we have chosen $\tau_{\mathrm{refine}} = 0.75$ and $\tau_{\mathrm{coarsen}} = \tau_{\mathrm{refine}}/3 = 0.25$.
In this simulation we selected $\gamma=4$ and have not employed the local substepping described in Section~\ref{sec:s-lts} (Equation~\eqref{eq:ratti-veroni-temporal}), as in this case its use did not result in an additional improvement of the solution quality. 
Note that the parameters selected for this benchmark are not the result of an extensive optimization procedure. \\ 
The results are reported in Fig.~\ref{fig:niederer-benchmark}, where we show that the proposed method leads to an accurate solution (the difference from the reference solutions is barely noticeable in the reported plot) while significantly reducing the computational load by a factor of about $12$.
This resulted in a wall-clock time of $\SI{105}{\second}$ while the reference solution using the Lie-Trotter-Godunov operator splitting scheme~\cite{QuGar:1999:aas,KriSarKlu:2013:nqo} required \SI{1237}{\second} ($\approx \SI{21}{\minute}$).

To investigate if the presented time step length adaptivity provides an advantage over combining a classical operator splitting scheme with spatial adaptivity, we have combined the spatial adaptivity presented in section~\ref{sec:spatial-adaptivity} with the Lie-Trotter-Godunov operator splitting. 
For simplicity we will call this scheme Lie-Trotter-Godunov+Adaptive Mesh Refinment (LTG+AMR) and used (conforming) finite elements in space. 
For the simulation of the conduction velocity benchmark~\cite{NieKerBenBerBerBraCheClaFenGarHeiLanMalPatPlaRodRoySacSeeSkaSmi:2011:vct} using the LTG+AMR scheme, we fixed $\Delta t = \SI{0.01}{\milli\second}$ and applied the refinement and coarsening steps only every \SI{0.15}{\milli\second}, in accordance with the setup for our S-LTS scheme, to keep the comparison fair. In addition we set the thresholds $\tau_{\mathrm{refine}} = 0.0025$ and $\tau_{\mathrm{coarsen}} = \tau_{\mathrm{refine}}/3 = 0.00083$ to control the mesh adaptation process in the LTG-AMR simulation.
Lower thresholds with respect to the ones chosen in the S-LTS scheme are selected to trigger mesh refinement in the LTG-AMR scheme earlier, i.e., for smaller errors. We note that for the same mesh size, the LTG-AMR scheme has fewer DOFs than the S-LTS DG scheme and therefore a lower threshold for mesh refinement in the LTG-AMR scheme (lowering $h_F$) leads to similar DOFs as in the S-LTS scheme.
Similar DOFs in both schemes lead to a similar accuracy (as suggested by numerical tests using the conduction velocity benchmark~\cite{NieKerBenBerBerBraCheClaFenGarHeiLanMalPatPlaRodRoySacSeeSkaSmi:2011:vct}) and therefore this strategy allows a fair comparison between schemes at equivalent accuracy.\footnote{In the LTG-AMR scheme, despite the change in spectrum due to AMR, the number of iterations in the inner linear solver during time marching did not increase significantly with respect to a solution computed without AMR and it remains between 5 and 9.}

The total wall-clock time for LTG-AMR scheme was $\SI{467}{\second}$, which is about $3$ times smaller than the wall-clock time required for the reference solution and 4 times larger than the time required by the newly proposed S-LTS scheme.
Since further optimization of each scheme parameters and implementation can affect the exact wall-clock time, the results presented here should be interpreted in terms of trends rather than exact speedups. 

\subsection{Spiral Wave Benchmark}
\label{sec:spiral-wave-benchmark}

The conduction velocity benchmark~\cite{NieKerBenBerBerBraCheClaFenGarHeiLanMalPatPlaRodRoySacSeeSkaSmi:2011:vct} primarely focuses on the convergence of the wave speed.
However, other features  of the computed solution may still be resolved incorrectly.
For example, despite obtaining a reasonable wave speed, it is possible that the wavefront is highly oscillatory or that the wave `attaches' to the grid, a phenomenon called lattice pinning~\cite[Fig.~39]{FenCheHasEva:2002:mms}. 
These incorrectly resolved features may affect the quality of simulated electrocardiograms or even result in a different qualitative behavior in re-entrant wave problems.
Formal error estimates like the Bochner norm, as used in our previous work~\cite{OgiPerBal:2023:sea}, are helpful to investigate detailed differences in the solutions.
However, this type of error estimates requires densely sampled time trajectories, easily leading to multiple terabytes of data for realistic problems and time scales, rendering this approach problematic.
Another downside of using error estimates based on the Bochner norm is that even small differences in the wave speed may lead to large errors while the overall wave propagation is preserved.
Due to these reasons we fall back on the standard visual comparison at a fixed time point, as often done in the community~(e.g.,~\cite{BelForBou:2014:pam,BelBriFor:2021:apa,MouPue:2021:dae,WooCanKal:2022:bcv}).
A visual comparison is consistent with the main goal of the benchmark proposed in this section, i.e., to show that the proposed method does not suffer from lattice pinning (see, e.g.,~\cite[Fig.~39]{FenCheHasEva:2002:mms}). %

\begin{figure}
    \centering
    \includegraphics[width=\textwidth]{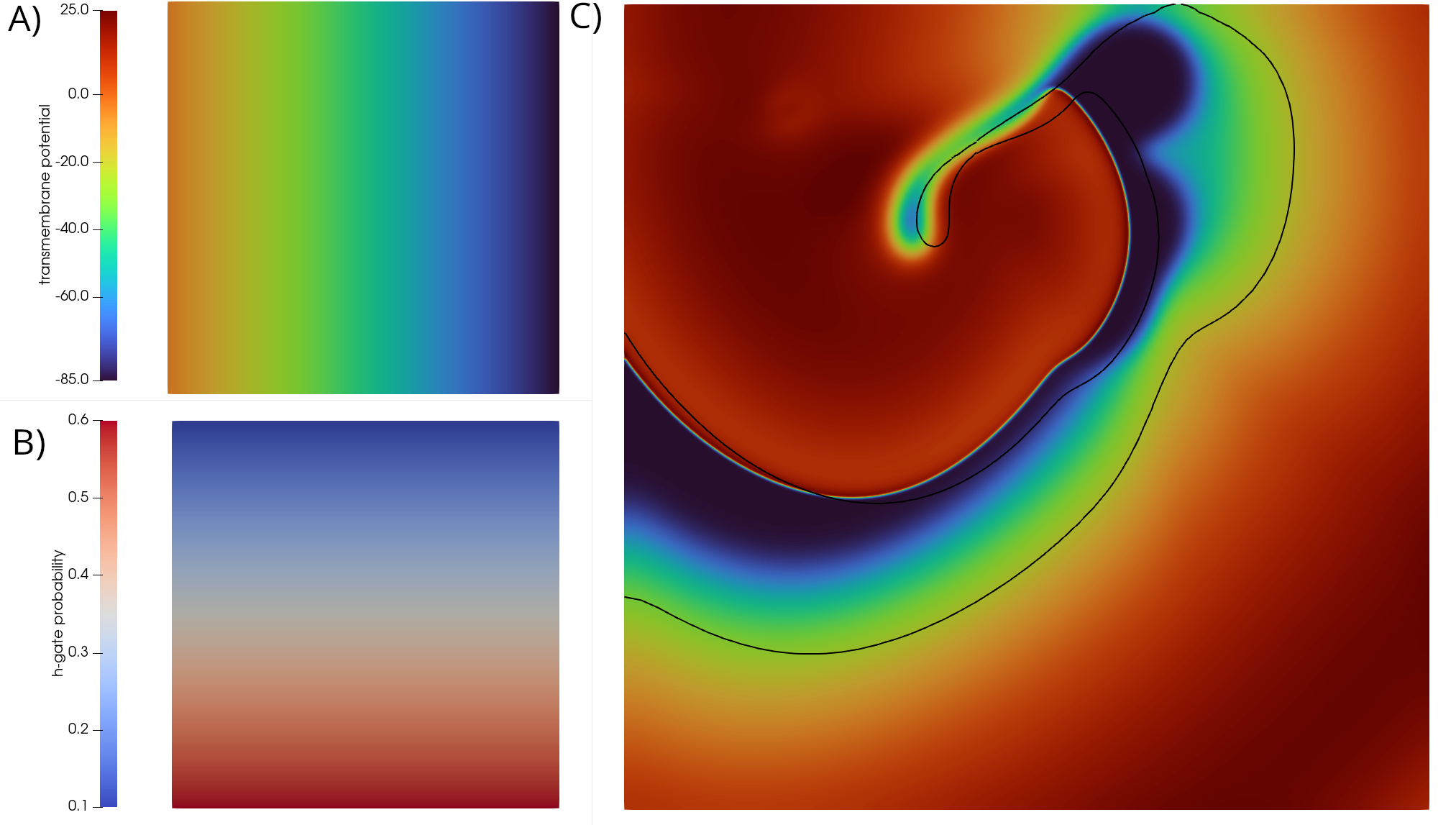}
    \caption{Spiral wave benchmark on a \SI{16}{\centi\meter} by \SI{16}{\centi\meter} domain. 
    A) Initial condition for the transmembrane potential.
    B) Initial condition for the h-gate.
    This setup is known to induce spiral waves across a large number of ionic models, (see, e.g.,~\cite{PatGalCorKabFenGra:2020:duq}).
    C) Transmembrane potential at $t=\SI{1}{\second}$ for the Lie-Trotter-Godunov operator splitting of the monodomain model with time step length $\Delta t=\SI{0.01}{\milli\second}$, grid size $\SI{0.156}{\milli\meter}$ and order $1$ Lagrange ansatz space. The legend to interpret the magnitude of the transmembrane potential is the same as the one reported in panel A.
    The black line represents the transmembrane isopotential line at $\SI{-30}{\milli\volt}$ for the newly proposed scheme. For the reference simulation and our scheme, we can observe that the wavefront shape is similar and the location of the wavefronts is close. Furthermore, no lattice pinning is observed.
    }
    \label{fig:spiral-wave-benchmark}
\end{figure}

In the benchmark presented here we use a quadratic domain of size \SI{16}{\centi\meter} $\times$ \SI{16}{\centi\meter} (spanned by the $x$- and $y$-coordinates), a simulation time equal to %
\SI{1000}{\milli\second}, and an isotropic diffusion tensor $\boldsymbol{D}=\mbox{diag}[0.1,0.1,0.1]$ $\SI{}{\milli\meter\squared\per\milli\second}$.
For the initial condition, we vary linearly $\varphi$ in the $x$ direction from \SI{10}{\milli\volt} to \SI{-85}{\milli\volt} and the h-gate state from 0.6 to 0.1 in $y$ direction.
The remaining ion channel states are set to their equilibrium state. 
This graded initial condition is known to likely induce a spiral wave~(e.g.,~\cite{PatGalCorKabFenGra:2020:duq}), as is the case in our numerical experiment (Fig.~\ref{fig:spiral-wave-benchmark}).
For the proposed S-LTS scheme, we use an initial coarse mesh with element size of \SI{1}{\centi\meter}, while the mesh for the reference simulation is based on an element size of \SI{0.156}{\milli\meter}.
In the S-LTS scheme, we further set the penalty $\gamma=8$ and order $2$ Lagrange polynomials for the ansatz space.
To control the scheme spatial adaptivity we have fixed $\tau_{\textrm{refine}}=1.0$ and $\tau_{\textrm{coarsen}} = \tau_{\textrm{refine}}/3 = 0.33$, while we set $\tau_{\mathrm{cell}} = 0.5$ and a substepping time step length of $\overline{\Delta t} = \SI{0.01}{\milli\second}$ to control the scheme temporal adaptivity.
In the reference solution, we use linear Lagrange polynomials for the ansatz space as this choice resulted in a faster simulation.

The wavefront computed in the reference and S-LTS solutions at the end of the simulation is presented in Fig.~\ref{fig:spiral-wave-benchmark}, from which it can be seen that the proposed method follows closely the reference solution.
Small differences in the wavefronts can be expected on such long time scales (the spiral wave already completed a full revolution) due to small differences in the wave speed at each time step, together with the problem that the AMR does not resolve accurately the first few milliseconds of the simulation --- a problem that is not addressed by the presented error estimators --- resulting in an initial offset between the waves.
In this case, the reference simulation required \SI{95100}{\second} ($\approx$ \SI{26}{\hour}), while the newly proposed S-LTS scheme required \SI{4135}{\second} ($\approx$ \SI{1}{\hour}).
We note that the wall-clock time needed to complete the reference and S-LTS simulations will be closer if a first order Lagrange ansatz was adopted in the S-LTS scheme as this choice will likely trigger additional spatial and temporal refinement during the solution steps.
However, the purpose of this benchmark is not to directly compare wall-clock time performance, but merely to show that no lattice pinning occurs when the solution is computed with the newly proposed scheme in an easily reproducible setup.

\subsection{Idealized Left Ventricle Benchmark}
\label{sec:lv-benchmark}

In this additional benchmark based on an idealized left ventricle, we aim to show that the proposed framework could be employed on more complex geometries than the rectangular domains (in 2D and 3D) used in the previous tests. 
The domain selected for this last benchmark is a truncated ellipsoid centered at the coordinate $(0,0,0)$, where the size of the inner and outer ellipsoids are \SI{6.5}{\milli\meter} $\times$ \SI{16.25}{\milli\meter} and \SI{11.0}{\milli\meter} $\times$ \SI{19.25}{\milli\meter}, respectively.
The ellipsoid is truncated at $x_3 \geq \SI{1.0}{\milli\meter}$. 
This model roughly approximates the size of a rabbit left ventricle.
For the first experiments we used a semi-structured hexahedral discretization of the domain.
The diffusion tensor is set to $\boldsymbol{D}=\mbox{diag}[0.1334,0.1334,0.1334]$ $\SI{}{\milli\meter\squared\per\milli\second}$.
The wave is initiated by a spherical stimulus with radius~\SI{10.0}{\milli\meter} and centered at $(0,0,0)$. Starting from its center and along the radial direction, this stimulus decreases linearly from \SI{100.0}{\micro\ampere\per\cubic\milli\meter} to \SI{0.0}{\micro\ampere\per\cubic\milli\meter}. 
Furthermore, it is applied from $t=\SI{0}{\milli\second}$ and decreases linearly over time until $t=\SI{2.0}{\milli\second}$.
This models a smooth activation of a large basal region of the idealized ventricle and, although it does not represent a physiological activation, it provides a clear activation wavefront to investigate the performance of the proposed scheme. 
We consider a semi-structured hexahedral mesh and Lagrange ansatz functions of order $1$ for all simulations.
The mesh for the reference solution has a fixed element size between $\SI{0.15}{\milli\meter}$ and $\SI{0.3}{\milli\meter}$ and we used a time step length of $\SI{0.001}{\milli\second}$ to ensure high resolution of the wave speed.
This solution was used to compare activation times between the reference and the S-LTS solutions. In addition, we recomputed the reference solution using a larger time step length of $\SI{0.01}{\milli\second}$ and used the (shorter) reference solution time when assessing the performance against the newly proposed scheme (this penalizes the newly proposed scheme, but offers a more sensible comparison for practical problems).
For our novel scheme we use $\tau_{\mathrm{refine}} = 0.5$, $\tau_{\mathrm{coarsen}} = \tau_{\mathrm{refine}}/10 = 0.05$, and $\tau_{\mathrm{cell}} = 0.5$ to control the spatial and temporal adaptivity, and a barrier time step length $\Delta t=\SI{0.1}{\milli\second}$.
The root mesh for the STL scheme contains approximately 4000 elements with size between $\approx \SI{0.6}{\milli\meter}$ and $\SI{1.0}{\milli\meter}$.
We fix the penalty parameter $\gamma=4$ and note that choosing a higher penalty $\gamma=8$ did not result in significantly higher accuracy.
However $\gamma=8$ decreased the CFL bound (eq.~\ref{eq:cfl-monodomain}) and hence the performance of our scheme.
Conversely, choosing smaller penalties (e.g., <<4) leads to unstable simulations.

In this benchmark we compare both the required solution time and the electrical wave conduction velocity.
In the cardiac simulation community it is well known that an important source of error is the over- or underestimation of the electrical wave conduction velocity~(see, e.g., \cite{PezHakSun:2015:sea} for a detailed analysis).
The most common measure of the wavefront velocity is the so-called local activation time (LAT).
The LAT can be computed as the time at every point $\boldsymbol{x}$ corresponding to the first instance when the transmembrane potential exceeds a critical threshold of \SI{-30}{\milli\volt}:
\begin{equation}
    \textrm{LAT}(\boldsymbol{x}) := \min_{t}\set*{ t \in [0;T] \mid \hat{\tilde{\boldsymbol{\varphi}}}_{\mathrm{m}}(\boldsymbol{x},t) > \SI{-30}{\milli\volt} }.
\end{equation}
Here $\hat{\tilde{\boldsymbol{\varphi}}}_{\mathrm{m}}$ %
denotes the piecewise linear interpolations in time of all $\hat{\tilde{\boldsymbol{\varphi}}}_{\mathrm{m}}^{n}$. %

\begin{figure}
    \centering
        \includegraphics[width=0.98\textwidth]{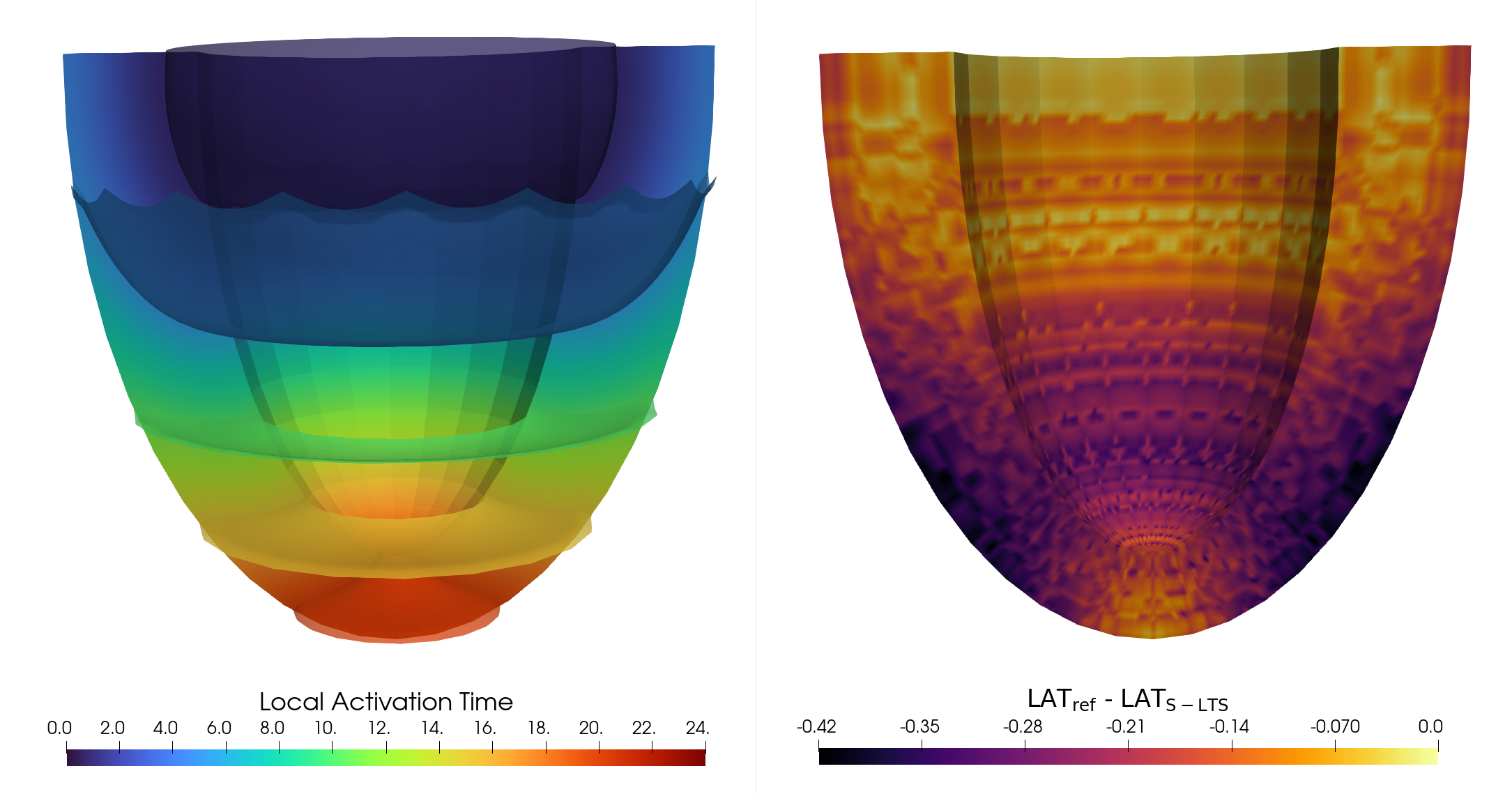}
    \caption{Benchmark on idealized left ventricle modeled as a truncated ellipsoid.
        \textbf{Left:} Reference activation time on a fine grid with element size between $\SI{0.15}{\milli\meter}$ and $\SI{0.3}{\milli\meter}$, and a small time step length of $\Delta t=\SI{0.001}{\milli\second}$ together with the isocontours at $t=1, 6, 12, 18, 24 \si{\milli\second}$.
        \textbf{Right:} Difference in activation times between the reference solution shown on the left and the newly proposed S-LTS scheme.
        The small artifacts results from a combination of a low sampling rate for the S-LTS scheme (the S-LTS solution is sampled only at the barrier time step of \SI{0.1}{\milli\second}), as well as a slightly varying wave speed in our scheme. For example, we notice a slight increase of the wavespeed for the S-LTS scheme toward the apical region. Nevertheless, the difference in activation times remains below $\SI{0.5}{\milli\second}$.
    }
    \label{fig:idealized-lv}
\end{figure}

The reference simulation required $\SI{3132}{\second}$ ($\approx \SI{50}{\minute}$) wall-clock time, while the newly proposed scheme required $\SI{1611}{\second}$ ($\approx \SI{25}{\minute}$), resulting in a speedup of only $2\times$.
The expected activation sequence and the difference in local activation timings is presented in Fig.~\ref{fig:idealized-lv}, showing that the activation time difference remained below 2\% at every point in the domain. 
We want to note that the speedup will significantly increase, if we also compute the repolarization times, because the S-LTS scheme can take larger time steps and have fewer elements than the reference solution.
By increasing the indicator bound to $\tau_{\mathrm{refine}} = 1.0$, $\tau_{\mathrm{coarsen}} = \tau_{\mathrm{refine}}/10 = 0.1$, we could make the S-LTS scheme a third faster (\SI{1025}{\second} or $\approx \SI{17}{\minute}$), while the activation error rises slightly above 2\% (the maximum difference was \SI{0.5}{\milli\second}).
Hence, we can gradually increase the indicator bound to decrease the wall-clock time until a desired work-precision tradeoff is reached.

All simulations presented up to this point were based on either structured or semi-structured grids.
However, in principle, the proposed scheme can handle unstructured grids as well.
To test the S-LTS performance in this scenario, we employed an unstructured hexahedral discretization of the idealized left ventricle used in the previous simulations.
The geometry has been discretized with Gmsh~\cite{GeuRem:2009:g3f} by hexahedralizing a tetrahedral mesh generated with the Frontal-Delaunay algorithm.
In this case, the total number of root elements is 2724.
The indicator bounds are set to $\tau_{\mathrm{refine}} = 2.5$ and $\tau_{\mathrm{coarsen}} = \tau_{\mathrm{refine}}/10 = 0.25$.
In this configuration, the activation time error was approximately 2\% (the maximum difference was \SI{0.5}{\milli\second}), similarly to the previous results obtained with a semi-structured hexahedral discretization of the domain.
The wall-clock time for this simulation was \SI{1263}{\second} or $\approx \SI{20}{\minute}$. 

\section{Discussion}
\label{sec:discussion}

We have presented a novel space-time adaptive explicit local time stepping scheme to solve monodomain problems in cardiac electrophysiology.
The local time stepping has been derived on the basis of a discontinuous Galerkin discretization of the monodomain problem.
It can be seen as a strategy to localize the Reaction-Tangent Controller scheme~\cite{OgiPerBal:2023:sea} on a per-element basis.
Although discontinuous Galerkin discretizations have more DOFs than their continuous counterparts and the time step length degrades quadratically with the element size, it is still reasonably efficient due to the small diffusivity of monodomain problems and the strict time step requirement to resolve the wavefront.
The utilization of an adaptive mesh refinement procedure allows to reduce the number of DOFs.
Analogue to the scheme by~\citet{QuGar:1999:aas}, the efficiency of our method stems from the significant reduction of necessary cell model evaluations, together with a smaller number of diffusion operator evaluations due to the reduction of the total number of time steps.

\subsection{Limitations}
\label{sec:limitations}

The most significant limitation of the proposed method is that the time step length decreases quadratically with element size.
This issue becomes even more important in the presence of small conductivities since they lead to sharper wavefronts that require a finer mesh to be properly resolved.
If it is important to resolve the slow wave with high accuracy, then our scheme's performance will severely degrade due to small time steps resulting from very fine grids in the slow regions.
In comparison, in similar circumstances, AMR techniques based on fully implicit discretizations will only lead to an increased number of DOFs while keeping the time step length constant. %

A limitation related to the presented performance studies consists in the missing investigation of optimization and parallelization strategies. Regarding the former, possible optimizations include implementing a caching infrastructure to recycle assembled face and element matrices on elements that are not coarsened or refined. %
While we still observe a measurable speed up in the presented benchmarks, in our implementation all matrices are reassembled and cached at the beginning of each time step. %
With respect to parallelization strategies, for example, it is possible to parallelize over the inner loops of Algorithm~\ref{alg:multiqueue-s-lts} in a shared memory setting.
However, we cannot expect optimal performance by doing so because the algorithm is memory bound and constructing optimal cache access patterns for local time stepping algorithms is not straightforward.
When a lower order ansatz is adopted in the problem discretization, the algorithm speed may be additionally degraded due to increased memory access (e.g.~\cite{TurKroBan:2016:wdp}).
On the other hand, higher order ansatzes pose a different challenge as they lead to larger element matrices and additional calculations.
For example, in our preliminary tests, ansatz spaces of order higher than $3$ harmed the performance of the scheme because the computational load per element increases without significant benefit in solution quality.
Although the investigation of optimial parallelization strategies combined with the order of the ansatz is left for future investigation, a possibility consists in employing a matrix-free algorithm implementation.
We refer to~\cite{TurKroBan:2016:wdp,KolFisMinDonBroDobWarTomSheAbdBarBeaCamChaDudKarKarKerLanMedMerObaPazRatSmiSpiSwiThoTomTom:2021:eed} for a more detailed discussion on this matter.
It is also possible to utilize a standard distributed memory parallelization strategy based on a dynamic domain decomposition, as provided in a user friendly setting in MFEM~\cite{AndAndBarBraCamCerDobDudFisKolPazStoTomAkkDahMedZam:2021:mmf}.
Although this type of parallelization will also be memory bound, the limitation might be less pronounced compared to the shared memory approach.
At this point, it remains unclear which specific load balancing technique should be applied as smaller elements may generate a large computational load due to additional time steps being computed and the possibly unbalanced load due to cell model evaluations (e.g., in the presence of branching in the cell model or simply different activation states).
A weighted space-filling curve with the number of substeps on an element as weight might be a good first proxy for the expected computational load. Due to the complexity and variety of possible solutions, the parallelization of the presented scheme requires future in depth investigation.

As a final note, we emphasize that it is not clear what ``sufficient accuracy'' in the context of cardiac electrophysiology simulations means, as the transient effect of under-resolved electrical fields on clinical markers has not been investigated in detail (see, e.g.,~\cite{PatCorGra:2019:cuq}).
As this requirement becomes more clear for specific electrophysiology applications, the parameters of the current scheme may be tuned to achieve the required accuracy.

\subsection{Application Areas}

Our method further expands the variety of numerical methods for solving monodomain problems, naturally raising the question about its possible application areas.
An obvious application where the presented scheme shines is the simulation of re-entrant wavefronts on realistic heart geometries -- as for example appearing in Wolff-Parkinson-White syndrome and during fibrillations -- because the localized wavefronts only cover a comparably small part of the domain.
These scenarios still require (at a minimum) a monodomain formulation to be properly resolved, in contrast to simplified formulations as, for example, the popular Reaction-Eikonal model~\cite{NeiCamPraNieBisVigPla:2017:ece}.
Our new strategy might also be helpful in electromechanical simulations with mechano-electrical feedback since it allows a tight integration of mechanical models while maintaining computational efficiency. %
A full derivation of the monodomain equation on a moving domain, as for example found in electromechanical simulations, is associated with a change in the diffusion tensor due to the structural rearrangement of the microstructure directions (i.e., the fiber, sheetlet, and normal directions). This would induce frequent reassembly of the matrices in the diffusion subproblem. %
Due to the coarser mesh away from the wavefront, AMR approaches lead to a smaller number of elements (compared to their counterparts based only on optimized operator splitting strategies) and therefore faster reassembly in electromechanical simulations.
Furthermore, by combining spatial (AMR) and temporal adaptivity, the S-LTS scheme retains the increased efficiency due to the time step adaptivity typical of operator splitting schemes.
Although further improvements are possible (e.g., based on two independent meshes for solving the electrical and mechanical components), the proposed S-LTS scheme is a promising approach for electromechanical simulations as it optimizes both in space and time when a single mesh is used.

\subsection{Concluding Remarks}

Our novel framework formalizes the intuitive idea that the computational load in cardiac electrophysiology monodomain simulations is primarily located at the wavefront.
The proposed scheme stems from the work on local time stepping in wave propagation problems from computational fluid dynamics~\cite{Kri:2010:elt}, where a related issue has been solved and explored.
Although there are several open challenges with the newly proposed framework, the presented benchmarks show that it is a promising candidate to speed up monodomain simulations.
We hope that this framework sheds light onto the possibility of further leveraging the localized nature of electrical waves on cardiac tissues and may foster numerical simulation strategies to incorporate the multiscale features in space and time to speed up simulations in cardiac electrophysiology.

%% file: figures/amr.tikz
\tikzset{
  treenode/.style = {align=center, inner sep=0pt, text centered,
    font=\sffamily},
  level0/.style = {treenode, circle, draw=black,
    fill=none, minimum size=3mm, very thick},
  level1/.style = {treenode, circle, draw=black,
    fill=none, minimum size=2mm},
  level2/.style = {treenode, circle, draw=black,
    fill=none, minimum size=2mm}
}

%{\hspace{15mm} Fine Grid}\\[4mm]
\hspace*{5mm}
\begin{tikzpicture}[scale=0.5,every node/.style={minimum size=1cm},level distance=10mm,sibling distance=10mm]
        % --- Root element ---
        \begin{scope}[ scale=0.65,
                yshift=-10,every node/.append style={
                yslant=0.5,xslant=-1},yslant=-cot(-35), xscale=sin(35)
                ]
            \fill[white,fill opacity=0.9] (0,0) rectangle (4,4);
            \draw[step=40mm, black] (0,0) grid (4,4); %defining grids
            %\draw[step=2.5mm, black!50,thin] (3,1) grid (4,2);  %Nested Grid
            \draw[black,very thick] (0,0) rectangle (4,4);%marking border 
        \end{scope} 
        \node [level0] at (1, -1) {};
        
        % --- First refinement ---
        \begin{scope}[ xshift=120,  scale=0.65,
                yshift=-10,every node/.append style={
                yslant=0.5,xslant=-1},yslant=-cot(-35), xscale=sin(35)
                ]
            \fill[white,fill opacity=0.9] (0,0) rectangle (4,4);
            \draw[step=20mm, black] (0,0) grid (4,4); %defining grids
            %\draw[step=2.5mm, black!50,thin] (3,1) grid (4,2);  %Nested Grid
            \draw[black,very thick] (0,0) rectangle (4,4);%marking border 
        \end{scope} 
        
        \node [level0] at (5, -1) {} 
            child{  
                node [level1] {}
            } child{  
                node [level1] {}
            } child {  
                node [level1] {}
            } child {  
                node [level1] {}
            };
        
        % --- Second refinement ---
        \begin{scope}[ xshift=240,  scale=0.65,
                yshift=-10,every node/.append style={
                yslant=0.5,xslant=-1},yslant=-cot(-35), xscale=sin(35)
                ]
            \fill[white,fill opacity=0.9] (0,0) rectangle (4,4);
            \draw[step=20mm, black] (0,0) grid (4,4); %defining grids
            \draw[step=10mm, black!50,thin] (0,0) grid (2,2);  %Nested Grid
            \draw[black,very thick] (0,0) rectangle (4,4);%marking border 
        \end{scope} 
        
        \node [level0] at (9.25, -1) {}
            child{  
                node [level1] {}
            } child{  
                node [level1] {}
            } child {  
                node [level1] {}
                child{  
                    node [level2] {}
                } child{  
                    node [level2] {}
                } child {  
                    node [level2] {}
                } child {  
                    node [level2] {}
                }
            } child {  
                node [level1] {}
            };
        
    % Separator line 
    \draw [dashed] (-1,-0.5) -- (11,-0.5);
\end{tikzpicture}